\documentclass{article}

\usepackage{paralist}
\usepackage{amssymb}
\usepackage{amsfonts}
\usepackage{amsmath}
\usepackage{graphicx}
\usepackage{latexsym}
\usepackage{epsf}
\usepackage{epsfig}

\newtheorem{theorem}{Theorem}

\newtheorem{definition}{Definition}

\newtheorem{lemma}{Lemma}

\newtheorem{proposition}{Proposition}
\newtheorem{remark}{Remark}

\title{Genus and braid index associated to sequences of renormalizable Lorenz maps}

\begin{document}

\maketitle
\centerline{\scshape Nuno Franco}
\medskip
{\footnotesize
 \centerline{CIMA-UE and Department of Mathematics, University of \'{E}vora}
   \centerline{Rua Rom\~{a}o Ramalho, 59, 7000-671 \'{E}vora, Portugal}
} 

\medskip

\centerline{\scshape and Lu\'
{i}s Silva }
\medskip
{\footnotesize
 \centerline{CIMA-UE and Department of Mathematics, University of \'{E}vora}
   \centerline{Rua Rom\~{a}o Ramalho, 59, 7000-671 \'{E}vora, Portugal}
} %

\bigskip


\begin{abstract}
We describe the Lorenz links generated by renormalizable Lorenz
maps with reducible kneading invariant
$(K_f^-,K_f^+)=(X,Y)*(S,W)$, in terms of the links corresponding
to each factor. This gives one new kind of operation that permits
us to generate new knots and links from old. Using this result we
obtain explicit formulas for the genus and the braid index of this
renormalizable Lorenz knots and links. Then we obtain explicit
formulas for sequences of these invariants, associated to
sequences of renormalizable Lorenz maps with kneading invariant
$(X,Y)*(S,W)^{*n}$, concluding that both grow exponentially. This
is specially relevant, since it is known that topological entropy
is constant on the archipelagoes of renormalization.

\end{abstract}

\section{Introduction} Let $\phi _t$ be a flow on $S^3$ with
countably many periodic orbits $(\tau _n)_{n=1}^{\infty}$. We can
look to each closed orbit as a knot in $S^3$. It was R. f.
Williams, in 1976, who first conjectured that non trivial knotting
occur in the Lorenz system (\cite{Wi1}). In 1983, Birmann and
Williams introduced the notion of template, in order to study the
knots and links (i.e. finite collections of knots, taking into
account the knotting between them) contained in the geometric
Lorenz attractor (\cite{BW}).

\begin{figure}[ht]\label{charts}\center
  \epsfig{file=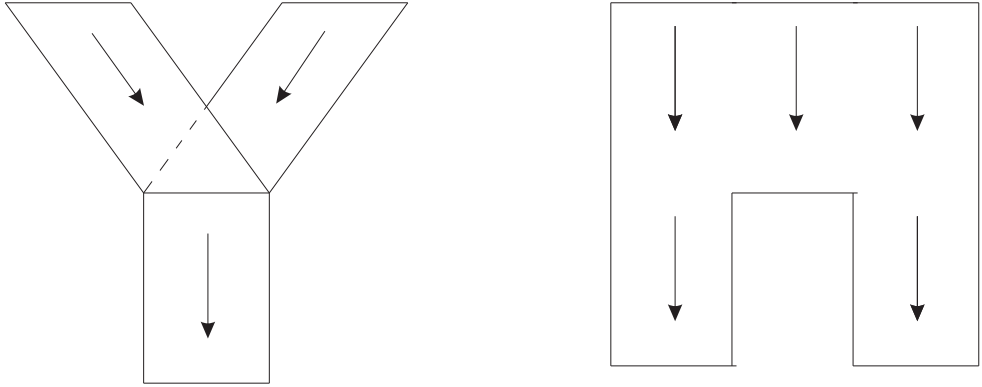,height=2in}\\
  \caption{Charts of templates: joining (left) and splitting (right)}
\end{figure}

A template, or knot holder, consists of a branched two manifold
with charts of two specific types, joining and splitting, together
with an expanding semiflow defined on it, see Figure \ref{charts}.
The relationship between templates and links of periodic orbits in
three dimensional flows is expressed in the following result,
known as Template Theorem, due to Birman and Williams in
\cite{BW}.

\begin{theorem}Given a flow $\phi _t$ on a
three-manifold $M$, having a hyperbolic chain-recurrent set, the
link of periodic orbits $L_{\phi}$ is in bijective correspondence
with the link of periodic orbits $L_{\mathcal{T}}$ on a particular
embedded template $\mathcal{T} \subset M$. On any finite sublink, this
correspondence is via ambient isotopy.
\end{theorem}

The dynamics of the semiflow on the Lorenz template are described
by the  first-return map to the branch line, which consists of a
one-dimensional map with one discontinuity, surjective and
strictly increasing  in each continuity subinterval.

What we now call a \textit{Lorenz flow} has a singularity of
saddle type with a one-dimensional unstable manifold and an
infinite set of hyperbolic periodic orbits, whose closure contains
the saddle point (see \cite{MeMA}). To describe the dynamics of
such a flow it is necessary to add a geometric hypotheses, just
like the one introduced in \cite{Wi} to study the original Lorenz
system. A Lorenz flow with this extra assumption is called a
\textit{Geometric Lorenz flow}. The dynamics of all these flows
can be described by  first-return one-dimensional maps with one
discontinuity, that are not necessarily surjective in the
continuity subintervals. This maps are called Lorenz maps, more
precisely, we will adopt the following definition introduced in
\cite{MeMA}.
\begin{definition}
Let $P<0<Q$ and $r \geq 1$. A $C^r$ Lorenz map $f:[P,Q]\rightarrow
[P,Q]$ is a map described by a pair $(f_-,f_+)$ where:
\begin{enumerate}
\item $f_-:[P,0]\rightarrow [P,Q]$ and $f_+:[0,Q]\rightarrow
[P,Q]$are continuous and strictly increasing maps;
\item $f(P)=P$, $f(Q)=Q$ and $f$ has no other fixed points in
$[P,Q]\backslash \{0\}$.
\item There exists $\rho >0$, the exponent of $f$, such that
$$f_-(x)=\widetilde{f}_-(|x|^{\rho}) \text{ and }
f_+(x)=\widetilde{f}_+(|x|^{\rho})$$ where  $\widetilde{f}_-$ and
$\widetilde{f}_+$, the coefficients of the Lorenz map, are $C^r$
diffeomorphisms defined on appropriate closed intervals.
\end{enumerate}
Because of the ambiguity at the point $0$, we consider the map
undefined in $0$. This Lorenz map is denoted by $(P,Q,f_-,f_+)$
(if there is no ambiguity about the interval of definition, we
erase the corresponding symbols $P,Q$).
\end{definition}

 In \cite{MeMA}, Martens and de Melo introduced some parametrized
 families of Lorenz maps that are universal in the sense that, given
 any geometric Lorenz flow, its dynamics is essentially the same
 as the dynamics of some element of the family, more precisely,
 consider  $\mathcal{L}^{r}$ the collection of all Lorenz maps of class
$C^r$. Endow $\mathcal{L}^r$ with a topology that takes care of the
domain, of the exponents and of the coefficients.

 \begin{definition} Let $\Lambda \subset \mathbb{R}^2$ be closed. A Lorenz
 family is a continuous map $F:\Lambda \rightarrow \mathcal{L} ^{r}$,
 $$F_{\lambda}=(P_{\lambda},Q_{\lambda},\varphi_{\lambda},\psi_{\lambda}).$$
\end{definition}

 A \textit{monotone Lorenz family} is a $C^3$ Lorenz family such
 that:
 \begin{enumerate}
 \item $F_{\lambda}$ has negative Schwarzian derivative for all
 $\lambda \in \Lambda$;
 \item $\Lambda =[0,1]\times [0,1]$;
 \item $F:(s,t)\rightarrow (-1,1,\varphi_s,\psi_t)$ and
 $\rho_{s,t}=\rho >1$;
 \item If $s_1<s_2$ then $\varphi_{s_1}(x)<\varphi_{s_2}(x)$ for
 all $x\in [-1,0]$ and if $t_1<t_2$ then $\psi_{t_1}(x)<\psi_{t_2}(x)$ for
 all $x\in [0,1]$;
 \item $\varphi_0(0)=0, \varphi_1(0)=1,\psi_0(0)=-1$ and
$ \psi_1(0)=0$;
\item $DF_{\lambda}(\pm 1)>1$ for all $\lambda \in \Lambda$.
\end{enumerate}

In \cite{MeMA} it is proved that Monotone Lorenz families are
full, in the sense that, if $F_{\lambda}$ is a monotone Lorenz
family, then for each given $C^2$ Lorenz map $f$ there is a
parameter $\lambda$ such that the dynamics of $F_{\lambda}$ are
essentially the same as the dynamics of $f$.

In \cite{H}, Holmes studied families of iterated horseshoe knots
which arise naturally associated to sequences of period-doubling
bifurcations of unimodal maps.

 It is well known, see for example \cite{DV}, that
period doubling bifurcations in the unimodal family are directly
related with the creation of a 2-renormalization interval, i.e. a
subinterval $J\subset I$ containing the critical point, such that
$f^2|_J$ is unimodal.

Basically there are two types of bifurcations in Lorenz maps (see
\cite{Pi}): the usual saddle-node or tangent bifurcations, when
the graph of $f^n$ is tangent to the diagonal $y=x$, and one
attractive and one repulsive $n$-periodic orbits are created or
destroyed; homoclinic bifurcations, when
$f^n(0^{\pm})=f^{n-1}(f_{\pm}(0))=0$ and one attractive
$n$-periodic orbit is created or destroyed in this way, these
bifurcations are directly related with homoclinic bifurcations  of
flows modelled by this kind of maps (see \cite{Pi}).

Considering a  monotone family of Lorenz maps, the homoclinic
bifurcations are realized in some lines in the parameters space,
called \textit{hom-lines} or \textit{bifurcation bones}.

It is known that (see \cite{Pi} and \cite{LSSR}), in the context
of Lorenz maps, renormalization intervals are created in each
intersection of two hom-lines. These points are called
\textit{homoclinic points} and are responsible for the
self-similar structure of the bifurcation skeleton of monotone
families of Lorenz maps. So it is reasonable to say that
\textit{homoclinic points are the
 Lorenz version of period-doubling bifurcation points}.

The idea of symbolic dynamics is to associate to each orbit of a
map, a symbolic sequence, called the \textit{itinerary} of the
corresponding point under the map. The pairs of sequences
corresponding to the orbits of the critical points determine all
the combinatorics of the map and are called \textit{kneading
invariants}; from this point of view, the renormalizability of a
map is equivalent to the reducibility of its kneading invariant as
the $\ast $-product of two other kneading invariants (see below
for the complete definitions).

At the topological and dynamical levels, we know a lot about  the
structure of renormalizable Lorenz maps, but, from the point of
view of knots and links generated by these maps, as far as we
know, this question was only superficially approached in \cite{W}.
So the objective of this work is to describe the structure and
invariants of knots and links generated by renormalizable Lorenz
maps with kneading invariants of type $(X,Y)*(S,T)$, by means of
the ones generated by $(X,Y)$ and $(S,T)$. Then we will study
sequences of invariants of knots and links generated by sequences
of kneading invariants corresponding to the iteration of the
$\ast$-product, i.e., corresponding to pairs of type
$$
(X,Y)\ast(K^{-},K^{+})^{\ast n}=(X(n),Y(n)),
$$
where $(K^{-},K^{+})^{\ast n}$ denotes the $\ast$-product of
$(K^{-},K^{+})$ with itself $n-1$ times.

Our main theorem describes the links corresponding to $n$-tuples
of (renormalizable) periodic itineraries of type
$((X,Y)*Z_1,\ldots ,(X,Y)*Z_n)$ in terms of the links
corresponding to $(X,Y)$ and $(Z_1, \ldots ,Z_n)$. This gives one
kind of operation that permits us to generate new knots and links
from old. Note that, unlike the case of Horseshoe knots, studied
by Holmes in \cite{H}, this operation do not corresponds to
cabling or any other operation that we know. Using this result we
proceed obtaining explicit formulas for the genus and the braid
index of knots and links corresponding to reducible sequences or
pairs of sequences depending, respectively, on the genus and on
the braid index of each factor. Then we obtain explicit formulas
for these invariants associated to pairs of type
$(X,Y)\ast(K^{-},K^{+})^{\ast n}=(X(n),Y(n))$, concluding that
both the genus and the braid index grow exponentially through
these sequences. This is specially relevant since it is known (see
\cite{LSSR}) that the topological entropy is constant in the
renormalization archipelagoes, and each of these sequences is
contained in one archipelago. So in this cases knot theory
provides much finer invariants for the classification of flows.

\section{Symbolic dynamics of Lorenz maps}

Denoting by $f^j=f\circ f^{j-1}$, $f^0=id$, the $j$-th iterate of
the map $f$, we define the \textit{itinerary} of a point $x$ under
a Lorenz map $f$ as $i_f(x)=(i_{f}(x))_j, j=0,1,\ldots$, where
\begin{equation*}
(i_{f}(x))_j=\left\{
\begin{matrix}
L & \text{if} & f^j(x)<0 \\
0 & \text{if} & f^j(x)=0 \\
R & \text{if} & f^j(x)>0
\end{matrix}
\right. .
\end{equation*}

It is obvious that the itinerary of a point $x$ will be a finite
sequence in the symbols $L$ and $R$ with $0$ as its last symbol,
if and only if $x$ is a pre-image of $0$ and otherwise it is one
infinite sequence in the symbols $L$ and $R$. So it is natural to
consider the symbolic space $\Sigma$ of sequences $X_{0} \cdots
X_{n}$ on the symbols $\{L,0,R\}$, such that $X_{i} \neq 0 $ for
all $i<n$ and: $n=\infty$ or $X_{n}=0$, with the lexicographic
order relation induced by $L<0<R$.

It is straightforward to verify that, for all $x,y \in [-1,1]$, we
have
\begin{enumerate}
\item If $x<y$ then $i_f(x)\leq i_f(y)$, and
\item If $i_f(x)<i_f(y)$ then $x<y$.
\end{enumerate}

We define the \textit{kneading invariant} associated to a Lorenz
map $f=(f_-,f_+)$, as
\begin{equation*}
K_{f}=(K_{f}^{-},K_{f}^{+})=(Li_f(f_-(0)),Ri_f(f_+(0))).
\end{equation*}%

We say that a pair $(X,Y)\in \Sigma \times \Sigma$ is \textit{admissible} if $%
(X,Y)=K_f$ for some Lorenz map $f$.

Consider the \textit{shift map} $s:\Sigma\setminus \{0\}
\rightarrow \Sigma$, $s(X_0\cdots X_n)=X_1\cdots X_n$. The set of
admissible pairs is characterized, combinatorially, in the
following way (see, for example, \cite{LSSR}).

\begin{proposition}
\label{t3} A pair $(X,Y)\in \Sigma \times \Sigma $ is admissible
if and only if $X_{0}=L$, $Y_{0}=R$ and, for $Z\in \{ X,Y \}$ we
have:
\newline (1) If $Z_{i}=L$ then $s ^{i}(Z)\leq X$;
\newline (2) If $Z_{i}=R$ then $s ^{i}(Z)\geq Y$;
with inequality (1) (resp. (2)) strict if $X$ (resp. $Y$) is
finite.
\end{proposition}

A sequence $X\in \Sigma $ is said to be $f-admissible$ if there
exists $x\in [ -1,1]$ such that $i_f(x)=X$. The $f$-admissible
sequences are completely determined by the kneading invariant
$K_f$ (see for example \cite{LSSR}), i. e., a sequence $X$ is
$f$-admissible if and only if it verifies the following
conditions:
\begin{enumerate}
\item If $X_i=L$ then $s^i(X)\leq K^-_f$;

\item If $X_i=R$ then $s^i(X)\geq K^+_f$;
with strict inequalities in the finite cases.
\end{enumerate}

\subsection{Renormalization and $*$-product}

In the context of Lorenz maps, we define renormalizability on the
following way, see for example \cite{Pi}:

\begin{definition} Let $f$ be a Lorenz map, then we say
that $f$ is $(n,m)$ renormalizable if there exist points
$P<y_{L}<0<y_{R}<Q$ such that
$$
g(x)=\left\{
\begin{array}{ll}
f^n(x) & \: if \: y_{L}\leq x<0 \\
f^m(x) & \: if \: 0<x\leq y_{R}%
\end{array}
\right.
$$
is  a Lorenz map.

The map $R_{(n,m)}(f)=g=(f^n,f^m)$ is called the
$(n,m)$-renormalization of $f$ and $[y_L,y_R]$ is the
corresponding renormalization interval.
\end{definition}

A sequence $X\in \Sigma$ is said to be \textit{maximal} if $X_0=L$
and $s^i(X)\leq X$ for all $i$ such that $X_i=L$, analogously a
sequence $Y\in \Sigma$ is \textit{minimal} if $Y_0=R$ and
$s^i(Y)\geq Y$ for all $i$ such that $Y_i=R$.

It is easy to prove that one infinite periodic sequence
$(X_0\cdots X_{m-1})^{\infty}$ with least period $m$ (the exponent
$\infty$ denotes the indefinite repetition of the sequence), is
maximal (resp. minimal) if and only if the finite sequence
$X_0\cdots X_{m-1} 0$ is maximal (resp. minimal).

Let $|X|$ be the length of a finite sequence $X=X_0\cdots
X_{|X|-1}0$, from the last observation it is reasonable to
identify each finite maximal or minimal sequence $X_0\cdots
X_{|X|-1}0$ with the corresponding infinite periodic sequence
$(X_0\cdots X_{|X|-1})^{\infty}$, this is the case, for example,
when we talk about the knot associated to a finite sequence.

It is also easy to prove that a pair of finite sequences
$$(X_0\ldots X_{|X|-1}0, Y_0\ldots Y_{|Y|-1}0)$$ is admissible, if
and only if the pair of infinite periodic sequences $$((X_0\cdots
X_{|X|-1})^{\infty},(Y_0\cdots Y_{|Y|-1})^{\infty})$$ is
admissible.

Considering a  monotone family of Lorenz maps, $F_{\lambda}$, the
homoclinic bifurcations are realized in the lines in the
parameters space such that the finite sequence $X$ is realized as
the left element of the kneading invariant, if $X$ is maximal  and
as the right element if $X$ is minimal. These lines are called
\textit{hom-lines} or \textit{bifurcation bones} and can be
defined as
$$
B(X)=\{\lambda\in \Lambda :K_{F_{\lambda}}^{-}=X\}.
$$
if $X$ is maximal and $$B(Y)=\{\lambda\in \Lambda
:K_{F_{\lambda}}^{+}=Y\}$$ if $Y$ is minimal.

\begin{figure}[ht]\label{bifdiagr}\center
  \epsfig{file=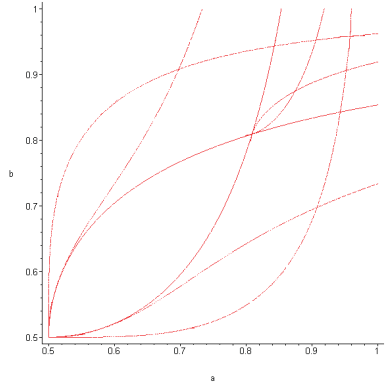,height=2in}\\
  \caption{Part of the bifurcation skeleton, namely, the bones with
   end point on the right side are the maximal bones corresponding to,
    from down to up, $LRL0$, $LR0$, $LRRL0=(LR0,RL0)*LR0$ and $LRR0$.
    With end point on the top we have the minimal bones corresponding,
    from left to right, to $RLL0$, $RL0$, $RLLR0=(LR0,RL0)*RL0$ and $RLR0$.
    The intersection of each two lines is the "vertex" of a similar copy of the all picture.}
\end{figure}

The union of the bifurcation bones is usually called the
\textit{bifurcation skeleton} (Figure \ref{bifdiagr}).

Obviously, two maximal or minimal bones corresponding to different
sequences can never intersect, so the only intersections we have
in the bifurcation skeleton are between maximal and minimal bones.
These points are called homoclinic and are where renormalization
intervals are created.

We define the $\ast $-product between a pair of
finite sequences $(X,Y)\in \Sigma \times \Sigma $, and a sequence
$U\in \Sigma $ as
$$
(X,Y)\ast U=\overline{U}_{0}\overline{U}_{1}\cdots
\overline{U}_{|U|-1}0,
$$
where
$$
\overline{U}_{i}=\left\{
\begin{array}{ll}
X_{0}\cdots X_{|X|-1} &\: if \: U_{i}=L \\
Y_{0}\cdots Y_{|Y|-1} &\: if \: U_{i}=R
\end{array}
\right. .
$$
Now we define the $\ast $-product between two pairs of sequences, $%
(X,Y),(U,T)\in \Sigma \times \Sigma $, $X$ and $Y$ finite, as
$$
(X,Y)\ast (U,T)=((X,Y)\ast U,(X,Y)\ast T).
$$

The next theorem states  that the reducibility relative to the
$\ast $-product is equivalent to the renormalizability of the map.
The proof can be found, for example, in  \cite{LSSR}.

\begin{theorem}
Let $f$ be a Lorenz map, then $f$ is $(n,m)$-renormalizable iff
there exist two admissible pairs $(X,Y)$ and $(U,T) $ such that
$|X|=n$, $|Y|=m$, $K_f=(X,Y)\ast (U,T)$ and
$K_{R_{(n,m)}(f)}=(U,T)$.
\end{theorem}

We also know from \cite{LSSR} that the product $(X,Y)\ast (U,T)$
is admissible if and only if both $(X,Y)$ and $(U,T)$ are
admissible, so for each finite admissible pair $(X,Y)$, the
subspace $(X,Y)*\{all\: admissible \: pairs\}$ is isomorphic to
the all space $\{all\: admissible \: pairs\}$, this provides a
self-similar structure in the symbolic space of kneading
invariants and, correspondingly, in the bifurcation skeleton. At
the topological and dynamical levels, the structure of the maps in
these similar subspaces is well described by renormalization and
$*$-product, but, from the point of view of knots and links
generated by these maps, as far as we know, this question was only
superficially approached in \cite{W}. So the objective of this
work is to describe the structure and invariants of knots and
links generated by Lorenz maps with kneading invariants of type
$(X,Y)*(S,T)$, relating it with the ones  generated by $(X,Y)$ and
$(S,T)$.

 First we will state some
useful properties of the $*$-product.

\begin{proposition}\label{l0}
Let $(X,Y)$ be one admissible pair of finite sequences, and $Z<Z'$,
then $(X,Y)*Z<(X,Y)*Z'$.
\end{proposition}

The proof is straightforward.

Denote by $X^{\infty}=(X_0\ldots
X_{|X|-1})^{\infty}=(X,Y)*L^{\infty}$ and by
$Y^{\infty}=(Y_0\ldots Y_{|Y|-1})^{\infty}=(X,Y)*R^{\infty}$. The
previous Proposition implies that $X^{\infty}\leq (X,Y)*Z \leq
Y^{\infty}$ for any sequence $Z$.

\begin{remark} From now on, we will freely identify the finite
sequence $X_0\ldots X_{|X|-1}0$, with the periodic sequence
$X^{\infty}=(X_0\ldots X_{|X|-1})^{\infty}$, wherever it is
convenient one or the other interpretation, for example, if
$p>|X|$ then we use $X_p$ to denote the element $X_{p\mod |X|}$.
\end{remark}

\begin{lemma}\label{l1}
Let $(X,Y)$ be one admissible pair of finite sequences, $0<q<|Y|$
and $Y_q=R$, then $Y_q\cdots Y_{|Y|-1}(X,Y)*Z \geq Y^{\infty}$,
for any sequence $Z$. Analogously, if $0<q<|X|$ and $X_q=L$, then
$X_q\cdots X_{|X|-1}(X,Y)*Z \leq X^{\infty}$, for any sequence
$Z$.
\end{lemma}
\textbf{Proof}

Since $(X,Y)$ is admissible, then $Y_q\ldots
Y_{|Y|-1}Y^{\infty}>Y^{\infty}$, so there exists $l$ such that
$Y_q\cdots Y_{q+l-1}=Y_0\ldots Y_{l-1}$ and $Y_{q+l}>Y_{l}$ . If
$q+l<|Y|$ the result follows immediately. If $q+l\geq|Y|$, then
necessarily $Y_{|Y|-q}=L$, because otherwise we would have
$Y^{\infty}>Y_{|Y|-q}\ldots$ and $Y_{|Y|-q}=R$, and this violates
admissibility. But then,
\begin{equation}\label{eq1}
Y_{|Y|-q}\cdots
Y_{|Y|-1}Y^{\infty}\leq
X^{\infty} \leq (X,Y)*Z
\end{equation}
and this gives the result. The proof of
the second part is analogous. $\blacksquare$

\begin{remark}Note that, equality in the first inequality of (\ref{eq1}), implies that $s^{|Y|-q}(Y^{\infty})=X^{\infty}$ so, if this is not the case, then the inequalities in the previous Lemma are strict.
\end{remark}

\begin{proposition}\label{l3}
Let $(X,Y)$ be one admissible pair of finite sequences and
$W,W'\in \{ X,Y \}$. If $s ^p(W^{\infty})< s ^q (W'^{\infty})$ and
$W_p \cdots W_{|W|-1} \neq W'_q \cdots W'_{|W'|-1}$ then
$$W_p \cdots W_{|W|-1}(X,Y)*Z \leq W'_q
\cdots W'_{|W'|-1}(X,Y)*Z'$$ for any sequences $Z,Z'$.
\end{proposition}

\textbf{Proof} The proof is divided in four cases:  $W=X$ and
$W'=Y$; $W=Y$ and $W'=X$; $W=W'=X$ and $W=W'=Y$. We will only
demonstrate specifically the first case, since the others follow
with analogous arguments..

Following the hypotheses, there exists $l$ such that $X_p\cdots
X_{p+l-1}=$ $ $ $Y_q \cdots Y_{q+l-1}$ and $X_{p+l}<Y_{q+l}$. If
$l<\min \{|X|-p,|Y|-q \}$, then the result follows immediately.

If $|X|-p\leq |Y|-q$ and $X_p\cdots X_{|X|-1}=Y_q\ldots
Y_{q+|X|-p-1}$, then $Y_{q+|X|-p}=R$, because otherwise we would
have $Y_{q+|X|-p}=L$ and $Y_{q+|X|-p}\ldots Y_{|Y|-1}Y^{\infty} >
X^{\infty}$, and this violates  admissibility of $(X,Y)$. So
$Y_{q+|X|-p}=R$ and, from Proposition \ref{l0} and Lemma \ref{l1},
\begin{equation}\label{eq2}
(X,Y)*Z\leq Y^{\infty} \leq
Y_{q+|X|-p}\ldots Y_{|Y|-1}(X,Y)*Z',
\end{equation} and the result follows.

If $|X|-p\geq |Y|-q$ and $X_p\cdots X_{p+|Y|-q-1}=Y_q\cdots
Y_{|Y|-1}$, then $X_{p+|Y|-q}=L$, because otherwise we would have
$X_{p+|Y|-q}=R$ and $X_{p+|Y|-q}\cdots < Y^{\infty}$, which
contradicts admissibility of $(X,Y)$. So $X_{p+|Y|-q}=L$ and
\begin{equation}\label{eq3}
X_{p+|Y|-q}\cdots X_{|X|-1}(X,Y)*Z\leq X^{\infty} \leq
(X,Y)*Z'\end{equation} and the result follows.

$\blacksquare$

\begin{remark} From the remark after Lemma \ref{l1}, we can only
have equalities in inequations \ref{eq2} and \ref{eq3} of the case
specifically studied
 and analogous in the other cases of the proof if $s^m(X^{\infty})=Y^{\infty}$
  for some $0<m <|X|$, so in the previous proposition, we can never have equality except if this happens.
\end{remark}

\begin{proposition}\label{l4}
Let $f$ be a $(n,m)$-renormalizable Lorenz map with
$R_{(n,m)}(f)=g$, renormalization interval $[y_L,y_R]$ and
kneading invariant $K_f=(X,Y)*(U,T)$, with $|X|=n$, $|Y|=m$ and
$K_g=(U,T)$. Then
$$i_f([y_L,y_R])=\{(X,Y)*Z\; such \; that \; Z \; is \;
g-admissible\}.$$
\end{proposition}

\textbf{Proof}

Note that $i_f(y_L)=(X_0\ldots
X_{|X|-1})^{\infty}=(X,Y)*L^{\infty}$ and $i_f(y_R)=(Y_0\ldots
Y_{|Y|-1})^{\infty}=(X,Y)*R^{\infty}$. Now, consider $x\in
[y_L,y_R]$, if $x<0$ then $(X,Y)*L^{\infty}\leq i_f(x) \leq
(X,Y)*U$, so the first $|X|$ symbols of $i_f(x)$ are equal to
$X_0\ldots X_{|X|-1}$. Analogously, if $y>0$ then the first $|Y|$
symbols of $i_f(x)$ are equal to $Y_0\ldots Y_{|Y|-1}$. Since
$f^n$ applies $[y_L,0[$ in to $[y_L,y_R]$ and $f^m$ applies
$]0,y_R]$ in to $[y_L,y_R]$, we can repeat the previous argument
to conclude that $i_f(x)=(X,Y)*Z$. The fact that $Z$ is
$g$-admissible and the reciprocal inclusion, follows immediately
from Proposition \ref{l0}.

$\blacksquare$

\section{Lorenz knots and links}

Let $n>0$ be an integer. We denote by $B_n$ the braid group on $n$
strings given by the following presentation:

$$B_n=\left\langle
\begin{array}{cc}
\sigma _{1},\sigma _{2},\ldots ,\sigma _{n-1} & \left|
\begin{array}{ll}
\sigma _{i}\sigma _{j}=\sigma _{j}\sigma _{i} & (|i-j|\geq 2) \\
\sigma _{i}\sigma _{i+1}\sigma _{i}=\sigma _{i+1}\sigma _{i}\sigma
_{i+1} & (i=1,\ldots ,n-2)
\end{array}
\right.
\end{array}
\right\rangle. \text{ }
$$
Where $\sigma_i$ denotes a crossing between the strings occupying
positions $i$ and $i+1$, such that the string in position $i$
crosses (in the up to down direction) over the other, analogously
$\sigma_i^{-1}$, the algebraic inverse of $\sigma_i$, denotes the
crossing between the same strings, but in the negative sense,
i.e., the string in position $i$ crosses under the other. A
\textit{positive braid} is a braid with only positive crossings. A
simple braid is a positive braid such that each two strings cross
each other at most once. So there is a canonical bijection between
the permutation group $\Sigma_n$ and  the set $S_n$, of simple
braids with $n$ strings, which associates to each permutation
$\pi$, the braid $b_{\pi}$, where each point $i$ is connected by a
straight line to $\pi (i)$, keeping all the crossings positive.

Let $X$ be a periodic sequence with least period $k$ and let
$\varphi\in \Sigma_k$ be the permutation that associates to each
$i$, the position occupied by $s^i(X)$ in the lexicographic
ordering of the $k$-tuple $(s(X),\ldots s^k(X))$ ($s^k(X)=X$).
Define $\pi \in \Sigma_k$ to be the permutation given by $\pi
(\varphi(i))=\varphi(i\mod k+1)$, i.e., $\pi(i)=\varphi(\varphi
^{-1}(i)+1)$. We associate to $\pi$ the corresponding simple braid
$b_{\pi}\in B_k$ and call it the \textit{Lorenz braid} associated
to $X$. Since $X$ is periodic, this braid represents a knot, and
we call it the \textit{Lorenz knot} associated to $X$.

\textbf{Example:} Let $X=(LRRLR)^{\infty}$. Hence we have
$s^5(X)=X$, $s(X)=(RRLRL)^{\infty}$, $s^2(X)=(RLRLR)^{\infty}$,
$s^3(X)=(LRLRR)^{\infty}$ and $s^4(X)=(RLRRL)^{\infty}$. Now after
lexicographic reordering the $s^i(X)$ we obtain $s^3(X)<
s^5(X)<s^2( X)<s^4( X)<s( X)$ and $\varphi=(1,5,2,3)$ written as a
disjoint cycle. Finally we obtain $\pi=(1,4,2,5,3)$ and
$b_{\pi}=\sigma_2\sigma_1\sigma_3\sigma_2\sigma_4\sigma_3$

\begin{figure}[ht]\center
  \epsfig{file=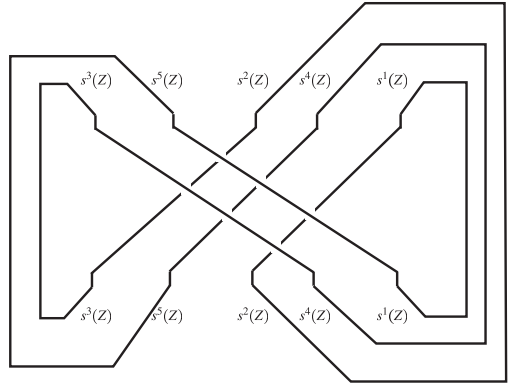,height=2in}\\
  \caption{The Lorenz knot associated to $X=(LRRLR)^{\infty}$}
\end{figure}

We can also generalize the previous algorithm to be used in the
case of a $p$-tuple of symbolic periodic sequences
$(X^1,\ldots,X^p)$ with periods $(k_1,\ldots,k_p)$. In this case
we proceed exactly as before with each one of the $X^j$. The
permutation $\varphi\in \Sigma_{k_1+\cdots+k_p}$ is the
permutation that describes the lexicographic ordering of the
$(k_1+\cdots k_p)$-tuple
$(s(X^1),\ldots,s^{k_1}(X^1),\ldots,s(X^p)\ldots,s^{k_p}(X^p))$
and $\pi\in \Sigma_{k_1+\cdots+k_p}$ is defined by
$\pi(\varphi(i))=\varphi(i+1)$ if there is no $q$ such that
$i=k_1+\cdots+k_q$ and $\pi(\varphi(i))=\varphi(k_1+\cdots
+k_{q-1}+1)$ if $i=k_1+\cdots k_q$, assuming $k_0=0$.

\begin{remark} What we are doing here is simply to mark in two parallel lines,
$k_1+\ldots +k_p$ points, corresponding in a ordered way, to the
sequences $s^{i_j}(X^j)$, $j=1,\ldots ,p$, $i_j=1,\ldots, k_j$ and
connect by straight lines the points corresponding to
$s^{i_j}(X^j)$ with the points corresponding to $s^{i_j+1}(X^j)$,
keeping the crossings positive.
\end{remark}

\section{The renormalization subtemplate}

A \textit{template} is a compact branched two-manifold with
boundary and a smooth expansive semiflow built locally from two
types of charts: joining and splitting (see Figure \ref{charts}).
Each chart carries a semiflow, endowing the template with an
expanding semiflow, and the gluing maps between charts must
reflect the semiflow and act linearly on the edges.

Following \cite{GHS}, we can take a semigroup structure on braided
templates. The generators of the braided template semigroup are:

\begin{enumerate}
\item $\sigma_i^{\pm}$, a positive (resp.negative) crossing between the strips
occupying the $i$-th and $(i+1)$-th positions;
 \item $\tau_i^{\pm}$, a half twist in the strip occupying the $i$-th position, in the positive (resp. negative) sense;
\item $\beta_i^{\pm}$, a branch line chart with the $i$-th and
$(i+1)$-th strips incoming , $2$ outgoing strips  and either a
positive ($\beta_i$) or negative ($\beta_i^-$) crossing at the
branch line.
\end{enumerate}

\begin{figure}[tbph]\label{gener}\center
  \epsfig{file=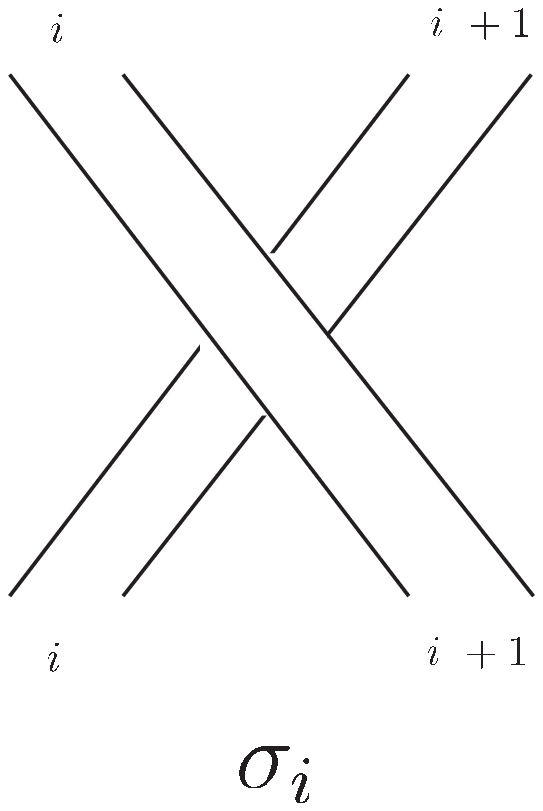,height=2in} \hspace{1.5cm}
  \epsfig{file=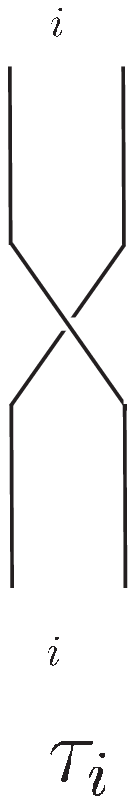,height=2in}\hspace{1.5cm}
  \epsfig{file=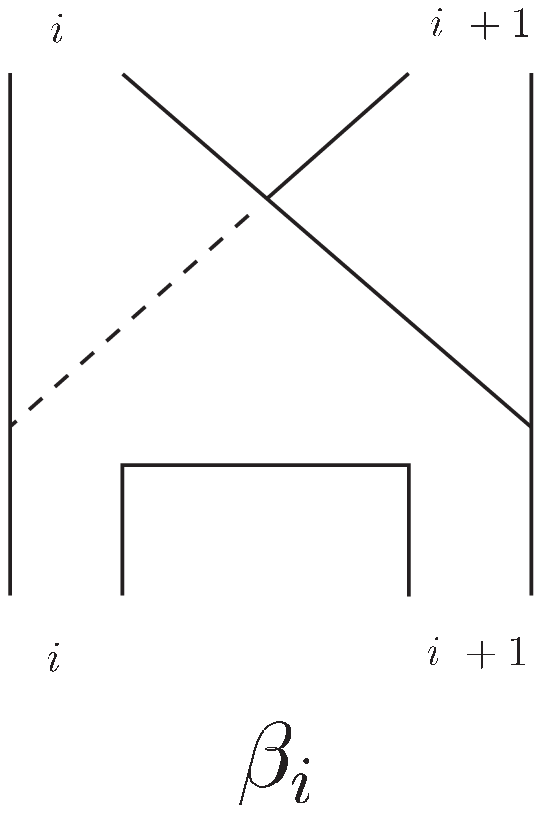,height=2in} \caption{Generators of
the braided template semigroup }
\end{figure}

 Given any pair of finite admissible sequences
$(X,Y)$, we define the \textit{tail's length} $m(X,Y)$ as
$$m(X,Y)=\min\{i\geq 0 : X_{|X|-1-i} \neq Y_{|Y|-1-i}\}$$

For a finite sequence $S$, let $n_L(S)=\# \{S_i : 0\leq i <|S|
\text{ and } S_i=L\}$, $n_R(S)=\# \{S_i : 0\leq i <|S| \text{ and
} S_i=R\}$.

  Now, to any finite admissible pair $(X,Y)$, we associate
a subtemplate $R (X,Y)$, the \textit{renormalization subtemplate}
associated to $(X,Y)$, on the following way: Consider the Lorenz
braid associated to $(X,Y)$, whose  word is $\sigma_{p_1}\cdots
\sigma_{p_k}$. Consider the relative position
$j=\varphi(|X|-m(X,Y))$, of $s^{|X|-m(X,Y)}(X)$.

If $s^n(X^{\infty})=Y^{\infty}$ for some $n<|X|$, then the
sequences $X$ and $Y$ generate the same Lorenz knot and we
consider the Lorenz link associated to $(X,Y)$ with only one
component. In this case $R(X,Y)$ is the template with $|X|$ strips
and word $\sigma_{p_1}\cdots \sigma_{p_k}\beta_j^{\pm}$, where the
signal $+$ in $\beta_j$ is taken if $X_{|X|-m(X,Y)-1}=L$  and the
signal $-$ is taken otherwise.

If $s^n(X^{\infty})\neq Y^{\infty}$ for all $n<|X|$ then $R(X,Y)$
is the template with $|X|+|Y|$ strips and word $\sigma_{p_1}\cdots
\sigma_{p_k}\beta_j^{\pm}$, where the signal $+$ in $\beta_j$ is
taken if $X_{|X|-m(X,Y)-1}=L$  and the signal $-$ is taken
otherwise.

\begin{figure}[tbph]\label{rentemp}\center
   \epsfig{file=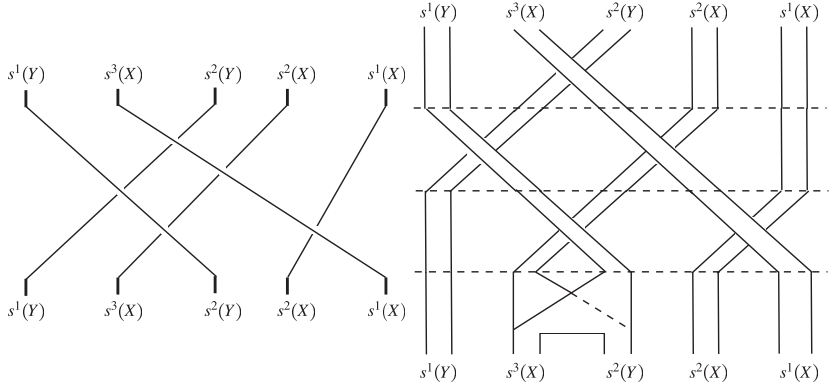,height=2in}\\
\caption{The renormalization template associated to $(X,Y)=((LRR)^{\infty},(RL)^{\infty}$}
\end{figure}

\begin{remark} What we are doing is simply to substitute each
string of the braid associated to $(X,Y)$ by a strip and add
$\beta_j^{\pm}$ according if $X_{|X|-m(X,Y)-1}=L$ or
$X_{|X|-m(X,Y)-1}=R$, respectively, see Figure 5
\end{remark}

The next theorem is naturally motivated from Proposition \ref{l4}.

\begin{theorem}
Let $(X,Y)$ be one admissible  pair of finite sequences and
$(Z^1,\ldots,Z^n)$ be a $n$-tuple of sequences whose associated
Lorenz link haves braid word $\sigma_{p_1}\cdots \sigma_{p_k}$,
then the Lorenz link associated to $((X,Y)*Z^1,\ldots,(X,Y)*Z^n)$
is the Lorenz link contained in $R(X,Y)$ with:
\begin{enumerate}
\item $|Z^1|+\cdots +|Z^n|$
strings in each strip if $s^i(X^{\infty})=Y^{\infty}$ for some
$i<|X|$.
\item  $n_L(Z^1)+ \cdots +n_L(Z^n)$ strings
in each strip associated to $X$ and $n_R(Z^1)+ \cdots +n_R(Z^n)$
strings in each strip associated to $Y$ if $s^i(X^{\infty})\neq
Y^{\infty}$ for all $i<|X|$.
\end{enumerate}
In both cases, the braid word of the restriction to the branch
line chart $\beta _j$ (respectively $\beta_j^-$) is
 $\sigma_{q+p_1}\cdots \sigma_{q+p_k}$ (respectively
$\sigma_{q+p_1}^{-1}\cdots \sigma_{q+p_k}^{-1}$), where $q+1$ is
the index of the left-most string getting in $\beta _j$.
\end{theorem}

\begin{figure}[tbph]\center
   \epsfig{file=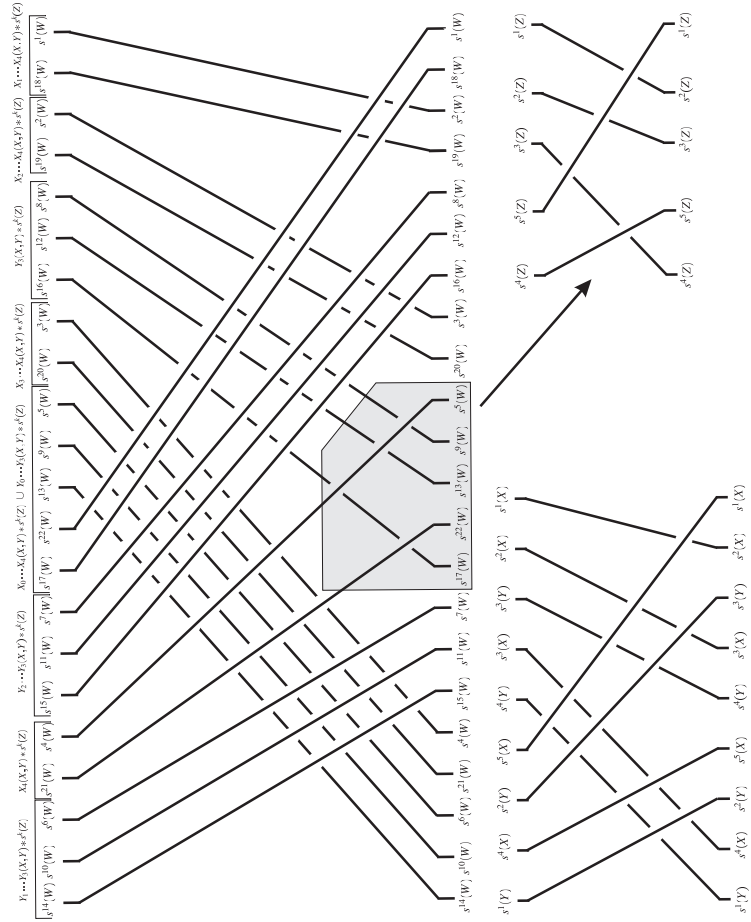, height=4.8in}\\
\caption{Pictoric ilustration of the Theorem: the Lorenz braid
associated to
$W=(X,Y)*Z=((LRRRL)^{\infty},(RLLR)^{\infty})*(LRRRL)^{\infty}$ on
the top, the Lorenz braids associated to $(X,Y)$ and $Z$ on the
bottom. }
\end{figure}

\newpage

\textbf{Proof}

We will only consider the case $s^i(X^{\infty})\neq Y^{\infty}$
for all $i<|X|$, since the proof of the other case is completely
analogous.

Without lost of generality, we can consider $n=1$, i.e. the Lorenz knot associated to $(X,Y)*Z$.

Consider the permutations $\varphi _{(X,Y)}$, $\varphi_{Z}$ and
$\varphi_{(X,Y)*Z}$ associated with the lexicographic ordering of
the sequences $(s(X),\ldots ,s^{|X|}(X),s(Y),\ldots ,s^{|Y|}(Y))$,
$(s(Z),\ldots ,s^{|Z|}(Z))$ and $(s((X,Y)*Z),\ldots ,
s^{|(X,Y)*Z|}((X,Y)*Z))$, respectively. Analogously consider $\pi
_{(X,Y)}$, $\pi_{Z}$ and $\pi_{(X,Y)*Z}$, the permutation induced
by the shift map over the respective lexicographically ordered
sequences  (see Section 3).

Let
$$W^k=\left\{\begin{array}{l}
  X \text{ if } Z_{\varphi_Z^{-1}(k)}=L\\
  Y \text{ if } Z_{\varphi_Z^{-1}(k)}=R
\end{array}\right.$$

For each $1\leq k \leq |Z|$ and $0\leq p \leq |W^k|-1$, define
$$\begin{array}{ll}
  \Phi(p,k) & =\varphi_{(X,Y)*Z}(s^p((X,Y)*s^{\varphi _Z^{-1}(k)}(Z))) \\
   & = \varphi_{(X,Y)*Z}(W_p^k\ldots W^k_{|W^k|-1}(X,Y)*s^{\varphi
   _Z^{-1}(k)+1}(Z)).
\end{array}
$$

From Propositions \ref{l0} and \ref{l3}, for each $1\leq i \leq
n_L(Z)-1$ and $0\leq p \leq |X|-m(X,Y)-1$, we have that
\begin{equation}\label{main1}
\Phi(p,i+1)=\Phi(p,i)+1, \end{equation}
 and, analogously, for each
$n_L(Z)+1\leq j \leq |Z|-1$ and $0\leq q \leq |Y|-m(X,Y)-1$
\begin{equation}\label{main2}\Phi(q,j+1)=\Phi(q,j)+1,\end{equation}
 This means that, in the lexicographic ordering
of $s^i((X,Y)*Z)$, the sequences $W_p\cdots W_{|W|-1}(X,Y)*s^k(Z)$
are all disposed together, constituting a set of $n_L(Z)$
sequences if $W=X$  and of $n_R(Z)$ sequences if $W=Y$, ordered by
the lexicographic ordering of $s^k(Z)$.

Moreover,
$$s(W_p\cdots W_{|W|-1}(X,Y)*s^{\varphi_Z^{-1}(k)+1}(Z))=W_{p+1}\cdots W_{|W|-1}(X,Y)*s^{\varphi_Z^{-1}(k)+1}(Z),$$ this means that
$$\pi_{(X,Y)*Z} (\Phi(p,k))=\Phi(p+1,k),$$
so there are exactly $n_L(Z)$ (resp. $n_R(Z)$) parallel strings from the set
$$\{\Phi(p,k),k=1,\cdots ,n_L(Z)\}\text{ to }\{\Phi(p+1,k),k=1,\cdots ,n_L(Z)\}$$
(resp. from
$$\{\Phi(p,k),k=n_L(Z)+1,\cdots ,|Z|\}\text{ to }\{\Phi(p+1,k),k=n_L(Z)+1,\cdots ,|Z|\}\text{)}.$$

Let us now consider the case $p\geq |W^k|-m(X,Y)$:

Since $\varphi_Z^{-1}(\pi_Z^{-1}(k))=\varphi_Z^{-1}(k)-1$, we have that
$$\pi_Z^{-1}(k)=\varphi_Z(\varphi_Z^{-1}(k)-1),$$
so, if $1\leq l\leq m(X,Y)$ then
$$\Phi(|W^{\pi_Z^{-1}(i)}|-l,\pi_Z^{-1}(i))=\varphi_{(X,Y)*Z}(W^{\pi_Z^{-1}(i)}_{|W^{\pi_Z^{-1}(i)}|-l}\ldots W^{\pi_Z^{-1}(i)}_{|W^{\pi_Z^{-1}(i)}|-1}(X,Y)*s^{\varphi_Z^{-1}(i)}(Z)).$$
From Propositions \ref{l0} and \ref{l3}, for each $1\leq i \leq
|Z|-1$ and $1\leq l \leq m(X,Y)$, we have that
\begin{equation}\label{main3}
\Phi(|W^{\pi_Z^{-1}(i+1)}|-l,\pi_Z^{-1}(i+1))=\Phi(|W^{\pi_Z^{-1}(i)}|-l,\pi_Z^{-1}(i))+1
\end{equation}
moreover, if $l>1$ then
$$\pi_{(X,Y)*Z}(\Phi(|W^i|-l,i)=\Phi(|W^i|-l+1,i).$$
So, in the lexicographic ordering of $s^i((X,Y)*Z)$, the sequences
$W_{|W|-l}\cdots W_{|W|-1}(X,Y)*s^k(Z)$ (with
$W=W^{\varphi_{Z}(k-1)}$) are all disposed together, constituting
a set of $|Z|$  sequences ordered according with $\pi_Z^{-1}(k)$
and there are exactly $|Z|$ parallel strings from the set
$$\{\Phi(|W^k|-l,k),k=1,\ldots ,|Z|\}\text{ to }\{\Phi(|W^k|-l+1,k),k=1,\ldots ,|Z|\}.$$ From Proposition \ref{l3}, the strings in the Lorenz link associated to $(X,Y)$, corresponding to
$s^{|X|-l}(X)\mapsto s^{|X|-l+1}(X)$ and $s^{|Y|-l}(Y)\mapsto
s^{|Y|-l+1}(Y)$ are parallel without any other string between
them. So we only have to divide each set of $|Z|$ strings from
$\{\Phi(|W^k|-l,k),k=1,\cdots ,|Z|\}$ to
$\{\Phi(|W^k|-l+1,k),k=1,\cdots ,|Z|\}$ in two subsets, the one on
the left with $n_L(X)$ strings contained in the strip associated
to $s^{|X|-l}(X)\mapsto s^{|X|-l+1}(X)$ and the one on the right
with $n_R(Z)$ strings contained in the strip associated to
$s^{|Y|-l}(Y)\mapsto s^{|Y|-l+1}(Y)$.

From Proposition \ref{l3}, $s^p(W^k)<s^q(W^{k'})$ implies that
$\Phi(p,k)<\Phi(q,k')$, so the sets $\{\Phi(p,k),k\}_p$ are
ordered according with $\varphi_{(X,Y)}$, this implies that the
transitions between this sets are done according with the
respective transitions in the Lorenz braid  associated to $(X,Y)$.
This finishes the proof of the first part of the theorem.

To see what happens in the branch line chart we must look to the transition to the tail, i.e. $l=m(X,Y)+1$:

If $X_{|X|-m(X,Y)-1}=L$ then $Y_{|Y|-m(X,Y)-1}=R$ and $\Phi(|X|-m(X,Y)-1,k)<\Phi(|Y|-m(X,Y)-1,k')$ for all $1\leq k\leq n_L(Z)$ and $n_L(Z)<k'\leq |Z|$.

On the other hand,
$$
s(W^k_{|W^k|-m(X,Y)-1}\ldots W^k_{|W^k|-1}(X,Y)*s^{\varphi_Z^{-1}(k)+1}(Z))=$$
\begin{equation}\label{eq4}\left\{\begin{array}{l}
     W^k_{|W^k|-m(X,Y)}\ldots W^k_{|W^k|-1}(X,Y)*s^{\varphi_Z^{-1}(k)+1}(Z) \text{ if } m(X,Y)>0 \\
       (X,Y)*s^{\varphi_Z^{-1}(k)+1}(Z) \text{ if }m(X,Y)=0
             \end{array}\right.\end{equation}
             and, while, from \ref{main1} and \ref{main2}, the elements $W^k_{|W^k|-m(X,Y)-1}\ldots W^k_{|W^k|-1}(X,Y)*s^{\varphi_Z^{-1}(k)+1}(Z)$
              are ordered according with $k$, from \ref{main3}, their shift images are ordered according with $\pi_Z^{-1}(k)$,
              this means that the permutation given by the strings connecting the two sets is exactly $\pi_Z$.

If $X_{|X|-m(X,Y)-1}=R$ then $Y_{|Y|-m(X,Y)-1}=L$ and $$\pi_{(X,Y)}(\varphi_{(X,Y)}(s^{|X|-m(X,Y)-1}(X))<\pi_{(X,Y)}(\varphi_{(X,Y)}(s^{|Y|-m(X,Y)-1}(Y)),$$
 this generates the crossing $\sigma_j$ in the braid associated to $(X,Y)$. Regarding to the braid associated to $(X,Y)*Z$, we have that
$$\Phi(|Y|-m(X,Y)-1,k')<\Phi(|X|-m(X,Y)-1,k)$$ for all $1\leq k\leq n_L(Z)$ and $n_L(Z)<k'\leq |Z|$. Now, considering
$$\xi (k)=\left\{\begin{array}{l}
                                                     n_L(Z)+k \text{ if } 1\leq k \leq n_R(Z) \\
                                                     k-n_R(Z) \text{ if } n_R(Z) < k \leq |Z|,
                                                   \end{array}\right.$$
while the elements $W^k_{|W^k|-m(X,Y)-1}\ldots
W^k_{|W^k|-1}(X,Y)*s^{\varphi_Z^{-1}(k)+1}(Z)$ are ordered
according with $\xi (k)$, their shift images are ordered according
with $\pi_Z^{-1}(k)$, this means that the permutation given by the
strings connecting the two sets, after crossing the $X$-strip with
the $Y$-strip, is exactly $\pi_Z$  so, from Remark \ref{trocas}
below, the braid restricted to the branch line chart is exactly
$\sigma _{p_1}^{-1} \cdots \sigma _{p_k}^{-1}$. $\blacksquare$

\begin{remark}\label{trocas}

The Reidemeister moves induce relations that are verified on the
braid group (resp. braided template semigroup). In the braided
template semigroup, one of these relations is
$\sigma_i\beta_i^-=\beta_i.$ This corresponds to make a
Reidemeister Type II move with the $i-$th and the $i+1-$th strips,
inverting the crossing on the $\beta_i$ line chart (changing the
sign), see Figure 7.

Given a simple braid $b$ (or the corresponding permutation) in a branch line chart $\beta_i$, it can be decomposed as the product of two braids. This decomposition $b=sb'$ is such that $s$ is a simple braid in the chart $\sigma_i$ and $b'$ is a mirrored simple braid (obtained from a simple braid changing all the crossings from positive to negative) in $\beta_i^-$.

\end{remark}

\begin{figure}[ht]\label{figlemma}\center
   \epsfig{file=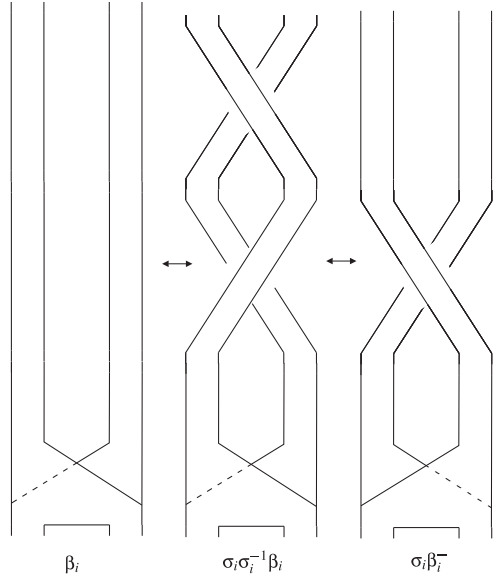, height=3in}\\
   \caption{Reidemeister Type II move and the relation $\sigma_i\beta_i^-=\beta_i.$}
\end{figure}

\section{Invariants}

We will start this section, introducing some terminology,
following \cite{BW}.

Let $\beta$ be a Lorenz braid:
\begin{enumerate}
\item The \textit{string index} is the number $n$ of strings in
$\beta$. It is the sum of the word lengths.
\item The \textit{braid index} of a knot is the minimum string index
among all closed braid representatives of that knot.
\item The \textit{crossing number} $c$ is the number of double
points in the projected image of the Lorenz braid $\beta$.
\item The \textit{linking number} $l(X,Y)$ is the number of crossings between one string from the knot associated to $X$
and one string from the knot associated to $Y$.
\item The \textit{genus} $g$ of a link $L$ is the genus of $M$,
where $M$ is an orientable surface of minimal genus spanned by
$L$.
\end{enumerate}

\begin{remark} Through all over this section we will only consider admissible pairs of finite sequences $(X,Y)$, such that $s^n(X^{\infty})\neq Y^{\infty}$ for all $n<|X|$. All the results respective to the case $s^n(X^{\infty})= Y^{\infty}$, follows analogously.
\end{remark}

\begin{lemma}\label{crossings}

Let $(X,Y)$ be a finite admissible pair and $S$ a finite sequence.
Then the crossing number $c((X,Y)*S)$ of the Lorenz braid
associated to $(X,Y)*S$, is given by
$$c((X,Y)* S)=c(X)n_L(S)^2+c(Y)n_R(S)^2+l(X,Y)n_L(S)n_R(S)\pm
c(S).$$ where we take the signal $+$ in $c(S)$ if
$X_{|X|-m(X,Y)-1}=L$ and
 the signal $-$ otherwise.

\end{lemma}

{\bf Proof: } There are four contributions to the computation of
$c((X,Y)* S)$. Two of them come from $c(X)$ and $c(Y)$, the third
one from $l(X,Y)$ and the fourth from $c(S)$. So $c(X)$ will be
counted $n_L(S)^2$ times and $c(Y)$ will be counted $n_R(S)^2$
times, this corresponds to the substitution of one crossing on $X$
(resp. $Y$) by the $n_L(S)^2$ (resp. $n_R(S)^2$) crossings arising
from inflating each $X$-string (resp. each $Y$-string) with
$n_L(S)$ (resp. $n_R(S)$)  strings. Similarly for each crossing
counted in $l(X,Y)$ we obtain $n_L(S)\times n_R(S)$ crossings.
Finally we must count the crossings in $\beta _j$, and, from the
Main Theorem and Remark \ref{trocas}, this means to  add or
subtract $c(S)$ according to if $X_{|X|-m(X,Y)-1}=L$ or
$X_{|X|-m(X,Y)-1}=R$. $\blacksquare$

\begin{lemma}\label{linking}

Let $(X,Y)$ and $(S,W)$ be  finite admissible pairs. Then the
linking number $l((X,Y)*(S,W))$ of the Lorenz braid associated to
$(X,Y)*(S,W)$, is given by

$$\begin{array}{l}
  l((X,Y)*(S,W))=l((X,Y)*S,(X,Y)*W)= \\
  2c(X)n_L(S)n_L(W)+2c(Y)n_R(S)n_R(W)+\\
  l(X,Y)(n_L(W)n_R(S)+n_R(W)n_L(S))\pm
l(S,W)
\end{array}$$

where we take the signal $+$ in $l(S,W)$ if $X_{|X|-m(X,Y)-1}=L$
and
 the signal $-$ otherwise.
\end{lemma}

{\bf Proof: } The proof is analogous to the proof of the previous
lemma , except that in this case we have to count the crossings in
the link $(X,Y)*(S,W)=((X,Y)*S,(X,Y)*W)$ between a string from
$(X,Y)*S$ and a string from $(X,Y)*W$. $\blacksquare$

Let $(S,W)$ be a finite admissible pair. We denote
$$A_{33}=\left[\begin{array}{ccc}
  n_L(S)^2 & n_L(W)^2 & 2 n_L(S)n_L(W) \\
  n_R(S)^2 & n_R(W)^2 & 2 n_R(S)n_R(W) \\
   n_L(S)n_R(S) &  n_L(W)n_R(W) &  n_L(W)n_R(S)+ n_L(S)n_R(W)
\end{array}\right]$$
 and $$ B_{13}=\left[\begin{array}{c}
  (n_L(S)+n_L(W))^2 \\
    (n_R(S)+n_R(W))^2 \\
    (n_L(S)+n_L(W))  (n_R(S)+n_R(W))
\end{array}\right]$$

\begin{lemma} \label{crossingslink}

Let $(X,Y)$ and $(S,W)$ be  finite admissible pairs. Then

$$c((X,Y)*(S,W))=\left[\begin{array}{ccc}
  c(X) & c(Y) & l(X,Y)
\end{array}\right]B_{13} \pm c((S,W))$$
where we take the signal $+$ in $c(S,W)$ if $X_{|X|-m(X,Y)-1}=L$
and
 the signal $-$ otherwise.
\end{lemma}

{\bf Proof: } The result follows from the previous two lemmas,
observing that
$c((X,Y)*(S,W))=c((X,Y)*S)+c((X,Y)*W)+l((X,Y)*(S,W))$. $\blacksquare$

\begin{remark}\label{erro}
We are concerned with the behavior of the invariants mentioned
above, through sequences of Lorenz braids and knots associated to
kneading invariants of type
$(A(n),B(n))=(X,Y)*(S,W)^n=(X,Y)*(S,W)^{n-1}*(S,W)$. Because of
the phenomenon described in Remark \ref{trocas} this invariants
may depend on the symbols $A(n-1)_{|A(n-1)|-m(A(n-1),B(n-1))-1}$,
so in the following we must consider four different cases:
\begin{enumerate}
\item If $X_{|X|-m(X,Y)-1}=S_{|S|-m(S,W)-1}=L$, then
$A(n)_{|A(n)|-m(A(n),B(n))-1}=L$ for all $n$.
\item If $X_{|X|-m(X,Y)-1}=L$ and $S_{|S|-m(S,W)-1}=R$, then
$$A(n)_{|A(n)|-m(A(n),B(n))-1}=\left\{\begin{array}{ll}
  L & \text{if } n \text{ is even} \\
  R & \text{if } n \text{ is odd}
\end{array}\right.$$.
\item If $X_{|X|-m(X,Y)-1}=R$ and $S_{|S|-m(S,W)-1}=L$, then
$A(n)_{|A(n)|-m(A(n),B(n))-1}=R$ for all $n$.
\item If $X_{|X|-m(X,Y)-1}=S_{|S|-m(S,W)-1}=R$, then
$$A(n)_{|A(n)|-m(A(n),B(n))-1}=\left\{\begin{array}{ll}
  L & \text{if } n \text{ is odd} \\
  R & \text{if } n \text{ is even}
\end{array}\right.$$.
\end{enumerate}
\end{remark}

\begin{lemma} \label{crossingslinkn}

Let $(X,Y)$ and $(S,W)$ be finite admissible pairs, then, for
$n\geq 2$ we have that

$$\begin{array}{l}
c((X,Y)*(S,W)^n)= \\
  \left( \left [\begin{array}{ccc} c(X) & c(Y)& l(X,Y) \end{array} \right] A_{33}^{n-1}
  + \left[ \begin{array}{ccc} c(S) & c(W)& l(S,W) \end{array} \right]
   \sum_{i=0}^{n-2} a_i A_{33}^i \right) B_{13} \\ + \alpha c(S,W)
   \end{array}$$
   where, considering the four cases of Remark \ref{erro}, we
   have: in Case 1 $a_i=1=\alpha$ for all $i$; in Case 2
   $a_i=(-1)^{i+n}$ and $\alpha = (-1)^{n+1}$; in Case 3 $a_i=-1=
   \alpha$ for all $i$; in Case 4 $a_i=(-1)^{i+n+1}$ and
   $\alpha=(-1)^n$.

\end{lemma}

{\bf Proof: } We will construct the formula recursively.

For $n=2$,  we want to compute
$c((X,Y)*(S,W)^2)=c((X,Y)*(S,W)*(S,W))$. Using first Lemma
\ref{crossingslink} and then Lemmas \ref{crossings} and
\ref{linking}, always considering Remark \ref{erro}, we have

$$\begin{array}{l}
  c((X,Y)*(S,W)*(S,W))=c(((X,Y)*S,(X,Y)*W)*(S,W))= \\
  \\
  \left[\begin{array}{ccc}
  c((X,Y)*S) & c((X,Y)*W) & l((X,Y)*(S,W)) \\
\end{array}\right]B_{13} \pm _1c(S,W)= \\
\\
  \left(\left [\begin{array}{ccc} c(X) & c(Y)& l(X,Y) \end{array}
\right] A_{33} \pm _0 \left[ \begin{array}{ccc} c(S) & c(W)&
l(S,W)
\end{array} \right]\right)B_{13} \pm _1 c(S,W)
\end{array}
$$
where $\pm _0$ depends on $X_{|X|-m(X,Y)-1}$ and $\pm _1$ depends
on $A(1)_{|A(1)|-m(A(1),B(1))-1}$.

For $n=3$,  we want to compute
$c((X,Y)*(S,W)^3)=c((X,Y)*(S,W)^2*(S,W))$. As in case $n=2$ we
obtain

$$\begin{array}{l}
  c((X,Y)*(S,W)^2*(S,W))= \\
  c(((X,Y)*(S,W)*S,(X,Y)*(S,W)*W)*(S,W))= \\
  \left[\begin{array}{ccc}
  c((X,Y)*(S,W)*S) & c((X,Y)*(S,W)*W) & l((X,Y)*(S,W)^{2}) \\
\end{array}\right]B_{13} \\
 \pm _2 c(S,W)= \\
\left [ \begin{array}{ccc} c((X,Y)*S) & c((X,Y)*W) &
l((X,Y)*(S,W))  \end{array} \right] A_{33}B_{13} \\
\pm _1 \left[\begin{array}{ccc} c(S) & c(W)& l(S,W) \end{array}
\right] B_{13} \pm _2 c(S,W)= \\
\left( \left [\begin{array}{ccc} c(X) & c(Y)& l(X,Y) \end{array}
\right] A_{33} \pm _0 \left[ \begin{array}{ccc} c(S) & c(W)&
l(S,W)
\end{array} \right] \right) A_{33}B_{13} \\
\pm _1 \left[  \begin{array}{ccc} c(S) & c(W)& l(S,W) \end{array}
\right] B_{13} \pm _2 c(S,W)=
 \\
 \left [\begin{array}{ccc} c(X) & c(Y)& l(X,Y) \end{array}
\right] A_{33}^2B_{13} \pm _0 \left[ \begin{array}{ccc} c(S) &
c(W)& l(S,W)
\end{array} \right] A_{33}B_{13}  \\
\pm _1 \left[  \begin{array}{ccc} c(S) & c(W)& l(S,W) \end{array}
\right] B_{13} \pm _2 c(S,W)
\end{array}
$$
where, considering $(A(0),B(0))=(X,Y)$, $\pm _i$ depends on
$A(i)_{|A(i)|-m(A(i),B(i))-1}$.

 We now obtain the formula recursively. $\blacksquare$

\begin{lemma} \label{length}

Let $n\in \mathbb{N}$ and $(X,Y)$ and $(S,W)$ be  finite admissible pairs. Then the string index of the Lorenz braid associated to $(X,Y)*(S,W)^n$ is

$$\left | (X,Y)*(S,W)^n \right |= \left[ |X||Y|\right]\left[ \begin{array}{cc}
  n_L(S) & n_L(W) \\
  n_R(S) & n_R(W)
\end{array} \right]^n \left[ \begin{array}{c}
  1 \\
  1
\end{array} \right]$$

\end{lemma}

{\bf Proof: } From the definition of $*$-product, we have that
$$\begin{array}{l}
  |(X,Y)*(S,W)|=|((X,Y)*S,(X,Y)*W)|=|(X,Y)*S|+|(X,Y)*W|= \\
  \\
  (n_L(S)+n_L(W))|X|+(n_R(S)+n_R(W))|Y|= \\
  \\
  \left[  \begin{array}{cc} |X| & |Y| \end{array} \right]\left[
\begin{array}{cc}
  n_L(S) & n_L(W) \\
  n_R(S) & n_R(W)
\end{array} \right]\left[ \begin{array}{c}
  1 \\
  1
\end{array} \right].
\end{array}
$$

Now, by induction,
$$\begin{array}{l}
  |(X,Y)*(S,W)^{n+1}|=|((X,Y)*S,(X,Y)*W)*(S,W)^n|=\\
  \\
  \left[  \begin{array}{cc} |(X,Y)*S| & |(X,Y)*W| \end{array} \right]\left[
\begin{array}{cc}
  n_L(S) & n_L(W) \\
  n_R(S) & n_R(W)
\end{array} \right]^n\left[ \begin{array}{c}
  1 \\
  1
\end{array} \right]=\\
 \\
  \left[  \begin{array}{cc} |X| & |Y| \end{array} \right]\left[
\begin{array}{cc}
  n_L(S) & n_L(W) \\
  n_R(S) & n_R(W)
\end{array} \right]^{n+1}\left[ \begin{array}{c}
  1 \\
  1
\end{array} \right].
\end{array}
$$
$\blacksquare$

\subsection{Braid index}
The \textit{trip number}, $t$, of a finite sequence $X$, is the
number of syllables in $X$, a syllable being a maximal subword of
$X$, of the form $L^aR^b$.

 Birmann and Williams conjectured in \cite{BW} that, for the case of a
 Lorenz knot $\tau$, $b(\tau)=t(\tau)$, where $t(\tau)$ is the trip number
  of the finite sequence associated to $\tau$. In \cite{W},  following a
  result obtained by Franks and Williams in \cite{FW}, Waddington observed that this conjecture is
  true. So our computations will be done about $t$.

\begin{proposition}[\bf Trip number and Braid index]\label{trip} Let $(X,Y)$ be a finite admissible pair, and $S$ be a finite sequence, then we have:
\begin{enumerate}

\item If $X_{|X|-1}=Y_{|Y|-1}$, then
$$t((X,Y)*S)=n_L(S)t(X)+n_R(S)t(Y).$$

\item If $X_{|X|-1}\neq Y_{|Y|-1}$, then
$$t((X,Y)*S)=n_L(S)t(X)+n_R(S)t(Y)\pm t(S),$$
where we take the signal $+$ in $t(S)$  if $X_{|X|-1}=L$ and
 signal $-$ otherwise.
\end{enumerate}

\end{proposition}

{\bf Proof: } The trip number is equal to the number of strings
that travel fom the $L$-side to the $R$-side (or, equivalently,
from the $R$-side to the $L$-side). In Case 1, since the branch
line chart $\beta _j$ is located completely in the $L$-side if
$X_{|X|-m(X,Y)}=L$ or in the $R$-side if $X_{|X|-m(X,Y)}=R$, the
only contributions to the trip number come from the strips
relative to $X$ and relative to $Y$ that travel from the $L$-side
to the $R$-side. Since there are exactly $n_L(S)$ strings in each
$X$-strip and $n_R(S)$ in each $Y$-strip, we get the result. In
the second case we have $j=n_L(X)+n_L(Y)$, this means that $\beta
_j$ haves one incoming strip from the $L$-side, other from the
$R$-side and the outgoing strips are also one in the $L$-side and
other in the $R$-side, so, from the Main theorem and Remark
\ref{trocas}, the $\beta _j$ chart will contribute with $\pm t(S)$
strings from the $L$-side to the $R$-side. $\blacksquare$

\begin{proposition} Let $(X,Y)$ and $(S,W)$ be finite admissible pairs. Then, for each $n\in \mathbb{N}$, we have:

\begin{enumerate}

\item If $X_{|X|-1}=Y_{|Y|-1}$, then
$$\left[\begin{array}{c}
  t((X,Y)*(S,W)^{n-1}*S) \\
    t((X,Y)*(S,W)^{n-1}*W)
\end{array}\right]= \left[ \begin{array}{cc}
  n_L(S) & n_R(S) \\
  n_L(W) & n_R(W)
\end{array}\right]^n \left[ \begin{array}{c}
  t(X) \\
  t(Y)
\end{array}\right].$$

\item If $X_{|X|-1}\neq Y_{|Y|-1}$ and $S_{|S|-1}\neq W_{|W|-1}$, then
$$\begin{array}{ll}
\left[\begin{array}{c}
  t((X,Y)*(S,W)^{n-1}*S) \\
    t((X,Y)*(S,W)^{n-1}*W)
\end{array}\right]= &
  \left[ \begin{array}{cc}
  n_L(S) & n_R(S) \\
  n_L(W) & n_R(W)
\end{array}\right]^n \left[ \begin{array}{c}
  t(X) \\
  t(Y)
\end{array}\right] \\
& \\
 &
  +
\sum_{i=0}^{n-1}a_i \left[ \begin{array}{cc}
  n_L(S) & n_R(S) \\
  n_L(W) & n_R(W)
\end{array}\right]^i \left[ \begin{array}{c}
  t(S) \\
  t(W)
\end{array} \right] \end{array}  $$
where, considering the cases from Remark \ref{erro},: in Case 1
$a_i=1$ for all $i$; in Case 2 $a_i=(-1)^{i+n+1}$; in Case 3
$a_i=-1$ for all $i$; in Case 4 $a_i=(-1)^{i+n}$.

\item
If $X_{|X|-1}\neq Y_{|Y|-1}$ and $S_{|S|-1}= W_{|W|-1}$, then
$$
\begin{array}{ll}
\left[\begin{array}{c}
  t((X,Y)*(S,W)^{n-1}*S) \\
    t((X,Y)*(S,W)^{n-1}*W)
\end{array}\right]
 = \left[ \begin{array}{cc}
  n_L(S) & n_R(S) \\
  n_L(W) & n_R(W)
\end{array}\right]^{n-1}
 \left[ \begin{array}{c}
  t((X,Y)*S) \\
  t((X,Y)*W)
\end{array}\right]= \\
 \left[ \begin{array}{cc}
  n_L(S) & n_R(S) \\
  n_L(W) & n_R(W)
\end{array}\right]^{n-1}
\left( \left[ \begin{array}{cc}
  n_L(S) & n_R(S) \\
  n_L(W) & n_R(W)
\end{array}\right]
 \left[ \begin{array}{c}
  t(X) \\
  t(X)
\end{array}\right] \pm
\left[ \begin{array}{c}
  t(S) \\
  t(W)
\end{array}\right] \right)
\end{array}
$$
where we take the signal $+$ in the last summand if $X_{|X|-1}=L$
and the signal $-$ otherwise.
\end{enumerate}

\end{proposition}

{\bf Proof: } This proof will be done by induction on $n$. The
case $n=1$ is just Proposition \ref{trip}. Now to prove the
induction step.

\begin{enumerate}

\item Suppose that the formula in case 1 is true for $n$. So  we want to compute

$$\begin{array}{l} \left[\begin{array}{c}
  t((X,Y)*(S,W)^{n}*S) \\
    t((X,Y)*(S,W)^{n}*W)
\end{array}\right]= \\ =\left[\begin{array}{c}
  t(((X,Y)*(S,W)^{n-1}*S,(X,Y)*(S,W)^{n-1}*W)*S) \\
  t(((X,Y)*(S,W)^{n-1}*S,(X,Y)*(S,W)^{n-1}*W)*W)
\end{array}\right]\end{array}$$

Hence from Proposition \ref{trip}: \noindent
$$\left[\begin{array}{c}
  t((X,Y)*(S,W)^{n}*S) \\
    t((X,Y)*(S,W)^{n}*W)
\end{array}\right]=$$
$$=\left[\begin{array}{c}
  n_L(S)t(((X,Y)*(S,W)^{n-1}*S)+n_R(S)t((X,Y)*(S,W)^{n-1}*W) \\
    n_L(W)t(((X,Y)*(S,W)^{n-1}*S)+n_R(W)t((X,Y)*(S,W)^{n-1}*W)
\end{array}\right]=$$
$$=\left[\begin{array}{cc}
  n_L(S) & n_R(S) \\
    n_L(W) & n_R(W)
\end{array}\right]\left[\begin{array}{c}
  t(((X,Y)*(S,W)^{n-1}*S) \\
    t((X,Y)*(S,W)^{n-1}*W)
\end{array}\right]$$

and we can apply our hypothesis to the second factor.

\item To prove case 2 we will follow the same steps as in case 1 to obtain

\noindent $$\left[\begin{array}{c}
  t((X,Y)*(S,W)^{n}*S) \\
    t((X,Y)*(S,W)^{n}*W)
\end{array}\right]=$$
$$\left[\begin{array}{c}
  n_L(S)t(((X,Y)*(S,W)^{n-1}*S)+n_R(S)t((X,Y)*(S,W)^{n-1}*W)\pm t(S) \\
    n_L(W)t(((X,Y)*(S,W)^{n-1}*S)+n_R(W)t((X,Y)*(S,W)^{n-1}*W) \pm t(W)
\end{array}\right]=$$
$$\left[\begin{array}{cc}
  n_L(S) & n_R(S) \\
    n_L(W) & n_R(W)
\end{array}\right]\left[\begin{array}{c}
  t(((X,Y)*(S,W)^{n-1}*S) \\
    t((X,Y)*(S,W)^{n-1}*W)
\end{array}\right]\pm \left[\begin{array}{c}
  t(S) \\
    t(W)
\end{array}\right]$$
where the signal $\pm$ in the last summand depends on Remark
\ref{erro}. Once again we obtain the desired result by applying
our hypothesis to the second factor of the first part of the sum .
\item Since $S_{|S|-1}=W_{|W|-1}$ then
$A(1)_{|A(1)|-1}=B(1)_{|B(1)|-1}$, so, because
$(X,Y)*(S,W)^n=((X,Y)*(S,W))*(S,W)^{n-1}$ we can apply part 1 of
this proposition to $((X,Y)*(S,W))*(S,W)^{n-1}$ and then part 2 of
Proposition \ref{trip} to $(X,Y)*(S,W)$ $\blacksquare$
\end{enumerate}

\subsection{Genus}

From Theorem 1.1.18 of \cite{GHS}, given a link $K$ and a braid
representative $b_K$ of the link, we have
\begin{equation}\label{eq5}g(K)=\frac{C-N-u}{2}+1,\end{equation} where $C$ is the number of
crossings in $b_K$, $N$ the string index and $u$ the number of
link components. We want now to compute $g((X,Y)* S)$.

\begin{proposition}[\textbf{Genus for knots}]\label{genus} Let $(X,Y)$ be a finite admissible pair and $S$ be a finite sequence. Then the genus of the knot associated to $(X,Y)*S$ is given by:
$$g((X,Y)*S)=$$
$$\frac{c(X)n_L(S)^2+c(Y)n_R(S)^2+l(X,Y)n_L(S)n_R(S)-n_L(S)|X|-n_R(S)|Y|+1\pm
c(S)}{2},
$$
where we take the signal $+$ in $c(S)$ if $X_{|X|-m(X,Y)-1}=L$ and
 the signal $-$ otherwise.

\end{proposition}

{\bf Proof: }  First notice that, because $(X,Y)* S$ is a knot we
have $u=1$. Now the number of strings in $(X,Y)* S$ is equal to
$|(X,Y)* S)|=n_L(S)|X|+n_R(S)|Y|$. The value of  $c((X,Y)* S)$ is
given by Lemma \ref{crossings}. So
$$
  g((X,Y)* S)=\frac{c((X,Y)* S)-|(X,Y)*S)|+1}{2}=$$

  $$\frac{c(X)n_L(S)^2+c(Y)n_R(S)^2+l(X,Y)n_L(S)n_R(S)-n_L(S)|X|-n_R(S)|Y|+1\pm
c(S)}{2}
$$

$\blacksquare $

\begin{proposition}[\textbf{Genus for links}]\label{genus2} Let $(X,Y)$ and $(S,W)$ be finite admissible pairs. Then the genus of the Lorenz link associated to $(X,Y)*(S,W)$ is given by:

$$g((X,Y)* (S,W))=\frac{c(X)(n_L(S)^2+n_L(W)^2)+c(Y)(n_R(S)^2+n_R(W)^2)}{2}+$$
$$\frac{l(X,Y)(n_L(S)n_R(S)+n_L(W)n_R(W))\pm(c(W)+c(S))+l((X,Y)*(S,W))}{2}-$$
$$\frac{(n_L(S)+n_L(W))|X|+(n_R(S)+n_R(W))|Y|}{2}$$

where we take the signal $+$ in $c(S)$ if $X_{|X|-m(X,Y)}=L$ and
 the signal $-$ otherwise.

\end{proposition}

{\bf Proof: } The proof is analogous to the proof of the previous
Proposition. $\blacksquare $

\begin{proposition} Let $(X,Y)$ and $(S,W)$ be finite admissible pairs. Then, for each $n\in\mathbb{N}$, the genus of the Lorenz link associated to $(X,Y)*(S,W)^n$ is given by:

$g((X,Y)*(S,W)^n)=$

$$\frac{1}{2}\left( \begin{array}{c}
  \left( \left [c(X)c(Y)l(X,Y) \right] A_{33}^{n-1} +
\left[ c(S)c(W)l(S,W) \right] \sum_{i=0}^{n-2}a_iA_{33}^i \right)B_{13} \\
\\
  + \alpha c(S,W)
-\left[ |X||Y|\right]\left[ \begin{array}{cc}
  n_L(S) & n_L(W) \\
  n_R(S) & n_R(W)
\end{array} \right]^n \left[ \begin{array}{c}
  1 \\
  1
\end{array} \right]
\end{array}   \right),$$
 where, considering the four cases of Remark \ref{erro}, we
   have: in Case 1 $a_i=1=\alpha$ for all $i$; in Case 2
   $a_i=(-1)^{i+n}$ and $\alpha = (-1)^{n+1}$; in Case 3 $a_i=-1=
   \alpha$ for all $i$; in Case 4 $a_i=(-1)^{i+n+1}$ and
   $\alpha=(-1)^n$.

\end{proposition}

{\bf Proof: } It is immediate, applying  the formulas in Lemmas
\ref{crossingslinkn} and \ref{length} in Equation \ref{eq5}.
$\blacksquare$

\end{document}